\documentclass[11pt]{article}
\usepackage{graphicx}
\usepackage{amsmath}
\usepackage{amsfonts}
\usepackage{epsfig}
\usepackage{setspace}

\newcommand{\be}{\begin{equation}}
\newcommand{\ee}{\end{equation}}
\newcommand{\beqn}{\begin{eqnarray}}
\newcommand{\eeqn}{\end{eqnarray}}
\newcommand{\beqns}{\begin{eqnarray*}}
\newcommand{\eeqns}{\end{eqnarray*}}

\newcommand{\EE}{\ensuremath{{\mathbb E}}}
\newcommand{\II}{\ensuremath{{\mathbb I}}}
\newcommand{\PP}{\ensuremath{{\mathbb P}}}

\newcommand{\fr}[1]{(\ref{#1})}

\newcommand{\eps}{\varepsilon}

\newtheorem{lemma}{Lemma}
\newtheorem{theorem}{Theorem}

\newtheorem{remark}{Remark}

\begin{document}
\allowdisplaybreaks
\title{\Large{\bf Minimax adaptive wavelet estimator for the anisotropic functional deconvolution model with unknown kernel}}

\author{
\large{ Rida Benhaddou}  \footnote{E-mail address: Benhaddo@ohio.edu} \ and\   Qing Liu  \\ \\
Department of Mathematics, Ohio University, Athens, OH 45701} 
\date{}

\doublespacing
\maketitle
\begin{abstract}
In the present paper, we consider the estimation of a periodic two-dimensional function $f(\cdot,\cdot)$ based on observations from its noisy convolution, and convolution kernel  $g(\cdot,\cdot)$ unknown. We derive the minimax lower bounds for the mean squared error assuming that $f$ belongs to certain Besov space and the kernel function $g$  satisfies some smoothness properties. We construct an adaptive hard-thresholding wavelet estimator that is asymptotically near-optimal within a logarithmic factor in a wide range of Besov balls. The proposed estimation algorithm implements a truncation to estimate the wavelet coefficients that is intended to stabilize the inversion. A limited simulations study confirms theoretical claims of the paper.\\


 \end{abstract} 

\section{Introduction.}
We consider the estimation problem of an unknown function $f(\cdot,\cdot)$ from observations $y(\cdot,\cdot)$ contaminated by Gaussian white noise in the following convolution model:
\be  \label{conv}
y(t,u)=\int_{0}^{1}f(x,u)g(t-x,u)dx+\eps z^{(1)}(t,u),\
g^{\delta}(t,u)=g(t,u)+\delta z^{(2)}(t,u),  
\ee
where $(t, u)\in [0,1]^2$, $g$ is the unknown blurring function with observations $g^{\delta}$, and $z^{(1)}$ and $z^{(2)}$ are two independent two-dimensional Gaussian white noises with covariance function  
\be \label{covz}
E\left[  z^{(k)}(t_1,u_1)z^{(k)}(t_2,u_2)  \right]=\delta(u_1-u_2)\delta(t_1-t_2), \ k=1,2,
\ee
and $\delta(\cdot)$ is the Dirac delta function. The parameters $\eps$ and $\delta$ are positive and satisfy $\eps, \delta \rightarrow 0$ asymptotically.\\
The discrete version of model \fr{conv} when $y(u,t)$ and $g^{\delta}(t,u)$ are observed at $NM$ points $(t_i,u_l)$, $l=1,\dots, M, i=1,\dots,N$, is as follows
\be \label{disc}
y(t_i,u_l)=\int_{0}^{1}f(x,u_l)g(t_i-x,u_l)dx+\sigma_1 z^{(1)}_{li},\
g^{\delta}(t_i,u_l)=g(t_i,u_l)+\sigma_2 z^{(2)}_{li} , 
\ee
 where $\sigma_1$ and $\sigma_2$ are two positive constants independent of $N$ and $M$, $u_l=l/M$ and $t_i=i/N$. The quantities $z^{(k)}_{li}$, with $k=1,2$, are zero-mean  i.i.d. normal random variables with  $\mathbb{E}\left[z^{(k)}_{l_1i_1}z^{(k)}_{l_2i_2}\right]=\delta(l_1-l_2)\delta(i_1-i_2)$. In addition, $z^{(1)}_{li}$ and $z^{(2)}_{li}$ are independent of each other.

Deconvolution model has witnessed a considerable number of publications since late 1980s and Donoho (1995) was the first to devise a wavelet solution to the problem. The list also includes Abramovich and Silverman (1998), Pensky and Vidakovic (1999), Walter and Shen (1999),  Johnstone, Kerkyacharian, Picard and Raimondo (2004), Donoho and Raimondo (2004), and Pensky and Sapatinas (2009). Functional deconvolution  has been investigated in Benhaddou, Pensky and Picard (2013), where they considered model \fr{conv} with $\delta=0$, which corresponds to the case when the kernel $g$ is known.  This model is motivated by experiments in which one needs to recover a two-dimensional function using  observations  of its convolutions along profiles $x=x_l$. This situation occurs, for example, in seismic inversions (see Robinson~(1999)). Another publication that is worth mentioning is that of Benhaddou, Pensky and Rajapakshage~(2019) who investigated the anisotropic functional Laplace deconvolution where the function under consideration is not periodic. 

In the present setting, the convolution kernel is unknown, but observations are available. This problem is referred to as the blind deconvolution. Inverse problem  with unknown operators in its general aspect was studied by Hoffmann and Reiss~(2008), where they proposed two nonlinear methods of estimating a one dimensional function. Delattre, Hoffmann, Picard and Vareschi (2012) considered the blind deconvolution problem and applied the block singular value decomposition (SVD) procedure to construct a wavelet estimator for a one-dimensional function that belongs to Sobolev space. Benhaddou (2018a) considered the blind deconvolution model under fractional Gaussian noise, where the function under consideration is univariate and periodic. The common feature between Hoffmann and Reiss~(2008) and  Benhaddou (2018a) is that they implement  preliminary thresholding procedures that eliminate the estimated wavelet coefficients that are judged to be too large. 

The objective of the paper is to construct an adaptive hard-thresholding wavelet estimator for model \fr{conv}. We focus on the regular-smooth convolution. Functional Fourier coefficients $g_m^{\delta}(u)$ with absolute values that vanish or are very close to zero at high frequencies may cause unstable inversion. For this reason, a preliminary stabilizing  thresholding procedure is applied to the functional Fourier coefficients of the `data' $g^{\delta}(t, u)$  to estimate the wavelet coefficients, taking advantage of the flexibility of the Meyer wavelet basis in the Fourier domain. That is, we apply the Meyer wavelet transform in the Fourier domain, and for each resolution level $j$, we  truncate the estimated wavelet coefficients at values $g_m^{\delta}(u)$ that are zero or close to zero. We show that the proposed approach is asymptotically near-optimal over a wide range of Besov balls under the $L^2$-risk. In addition, we demonstrate that the convergence rates are expressed as the maxima between two terms, taking into account both the noise sources.  Similar behavior has been pointed out in Hoffmann and Reiss~(2008), Vareschi~(2015) and Benhaddou~(2018a, 2018b). It should be noted that with $\delta=0$, our convergence rates coincide with those in Benhaddou et al.~(2013), and with $\delta=0$ and $p=2$, our convergence rates match those in Benhaddou~(2017). Finally, with $\alpha=1$, and with $\alpha_{1l}=1$ and $M=1$, our rates are comparable to those in Benhaddou~(2018a) and Benhaddou~(2018b), respectively.

The rest of the paper is organized as follows. In section 2, we describe in details the estimation procedure. In section 3, we study the asymptotic performance of the proposed estimator in terms of its minimax squared loss. In Section 4, we consider a limited simulations study to assess the goodness of our estimator when the sample size is finite. Finally, Section 5 contains the proofs of our theoretical findings. 
\section{Estimation algorithm}
Let  $\langle \cdot,\cdot \rangle$ denote the inner product in the Hilbert space $L^2(U)$, where $U=[0,1]^2$, i.e., $\langle f,g \rangle= \int_{0}^{1}f(t,u)\overline{g(t,u)}dt$ for $f,g \in L^2(U)$, and $\overline{g} $ is the complex conjugate of $g$. Denote $e_m(t)=e^{i2\pi mt}$ as a Fourier basis on the interval $[0,1]$. Let  $y_m(u)=\langle e_m,y(\cdot,u) \rangle$, $z_m(u)=\langle e_m,z(\cdot,u) \rangle$, $g^\delta_m(u)=\langle e_m,g^\delta(\cdot,u) \rangle$, $g_m(u)=\langle e_m,g(\cdot,u) \rangle$, $f_m(u)=\langle e_m,f(\cdot,u) \rangle$ be functional Fourier coefficients of functions $y, z, g^\delta, g, f$, respectively. Applying Fourier transform to equation \fr{conv} we get
\be \label{fconv}
y_m(u)=f_m(u)g_m(u)+\eps z^{(1)}_m(u),
\
g^{\delta}_m(u)=g_m(u)+\delta z^{(2)}_m(u).
\ee
where $z^{(k)}_m(u)$ are generalized one-dimensional Gaussian processes such that
\be 
\EE\left[z^{(k)}_{m_1} (u_1)z^{(k)}_{m_2} (u_2)\right] = \delta_{m_1,m_2} \delta (u_1 - u_2),\   k=1, 2. \label{zm}
\ee
Consider a bounded bandwidth periodized wavelet basis (Meyer-type) and a finitely supported periodized $s_0$-regular wavelet basis (Daubechies-type). Denote the wavelet functions of these two bases by $\psi_{m_0-1,k}(t)$ and $\eta_{m'_0-1,k'}(u)$ respectively, where $m_0$ and $m'_0$ correspond to the lowest resolution levels for the two bases. Then, the function  $f$ has the wavelet series representation given by
\be \label{f}
f(t,u)=\sum\limits_{j=m_0-1}^\infty \sum\limits_{j'=m_{0'}-1}^\infty \sum\limits_{k=0}^{2^j-1} \sum\limits_{k'=0}^{2^{j'}-1} \beta_{j,k; j',k'}\psi_{j,k}(t)\eta_{j',k'}(u),
\ee 
with
\be \label{beta}
\beta_{j,k; j',k'}=\langle \langle f(t,u), \psi_{j,k} (t)  \rangle, \eta_{j',k'}(u) \rangle=\sum\limits_{m \in W_j} \int_{0}^{1} f_m(u)\ \eta_{j',k'}(u) du\ \overline{\psi}_{j,k,m} ,
\ee
and $\psi_{j,k,m} =\langle e_m(t), \psi_{j,k}(t)\rangle $ are the Fourier coefficients of $\psi_{j,k}$. It is well-known (see, e.g, Johnstone et al.~(2004), section 3.1) that under the Fourier domain and for any $j\geq j_0$, one has
\be
W_j=\left\{m: \psi_{j,k,m}\ne 0 \right\} \subseteq \frac{2\pi}{3} \left[-2^{j+2},-2^j\right] \cup \left[2^j,2^{j+2}\right], \label{wj}
\ee 
and the cardinality of $W_j$ is $|W_j|=4\pi \cdot 2^j$. If $g$ were known, the problem would reduce to the regular deconvolution model studied in Benhaddou et al.~(2013). However, since $g$ is unknown and contaminated with gaussian white noise, a preliminary thresholding procedure is to be applied to estimate the wavelet coefficients $\beta_{j,k; j',k'}$. Let us define the quantities
\be  \label{gdel}
\frac{1}{\hat{g}^{\delta}_m(u)}= 
\begin{cases}
\frac{1}{g^{\delta}_m(u)} & \text{if } \left|  g^{\delta}_m(u) \right|>\kappa \delta \sqrt{\ln(1/\delta) },
\\
0 & \text{if } \left|  g^{\delta}_m(u) \right| \leq \kappa \delta \sqrt{\ln(1/\delta) },
\end{cases}
\ee 
where $\kappa$ is a positive constant independent of $m$ and $\delta$. Then, using \fr{beta} and \fr{gdel}, and by Plancherel formula,  a truncated estimator for $\beta_{j,k; j',k'}$ is given by
\be  \label{bhat}
\hat{\beta}_{j,k;j',k'}=\sum\limits_{m \in W_j} \int_{0}^{1} \frac{y_m(u)}{\hat{g}^{\delta}_m(u)} \ \eta_{j',k'}(u) du\  \overline{\psi}_{j,k,m}.
\ee 
Bear in mind that thresholding \fr{gdel} enables us to have a stable inversion since it eliminates the values of $g^{\delta}_m(u)$ that are equal to zero or relatively close to zero. Now, define 
\be \label{ojj}
\resizebox{.93 \textwidth}{!}{$\Omega(J,J')=\{ \omega=(j,k; j',k'): m_0\leq j\leq J-1, m'_0\leq j'\leq J'-1, k=0,\cdots, 2^{j}-1,k'=0,\cdots, 2^{j'}-1 \}.$}
\ee
Then, consider the hard-thresholding estimator for $f(\cdot,\cdot)$ given by 
\be \label{fhat}
\hat{f}(t,u)=\sum\limits_{j=m_0-1}^{J-1} \sum\limits_{j'=m_{0'}-1}^{J'-1} \sum\limits_{k=0}^{2^j-1} \sum\limits_{k'=0}^{2^{j'}-1} \hat{\beta}_{\omega} \II\left(|\hat{\beta}_{\omega} |>\lambda^j_{\eps, \delta}\right)   \psi_{j,k}(t)\eta_{j',k'}(u),
\ee
where $J$ and $J'$ will be determined later, and $\lambda^j_{\eps, \delta}$ will be chosen based on moment properties of $\hat{\beta}_{\omega}$. Let us now introduce an assumption that pertains the smoothness property of the convolution kernel $g(t, u)$.\\
\noindent
{\bf Assumption A.1.  }\label{A1}Assume that the functional Fourier coefficients $g_m(u)$ of $g(t,u)$ are uniformly bounded from below and above, that is, there exists positive constants $\nu$, $C_1 $ and $C_2$, all independent of $m$ and $u$, such that 
\be 
C_1 |m|^{-2\nu} \leq |g_m(u)|^2 \leq C_2 |m|^{-2\nu}  .\label{gm}
\ee
The next step is to evaluate the mean-squared error of \fr{bhat}. Indeed, define for $j\geq m_0$ and positive constant $\rho$, the sets $\Omega_1$ and $\Omega_2$ such that
\beqn 
\Omega_1(j) &=&\left \{m\in W_j: \left|g^{\delta}_m(u) \right|>\kappa \delta \sqrt{\ln(1/\delta) } \right \}  ,  \label{o1}  \\
\Omega_2(j) &=&\left \{m\in W_j: \left| \delta z^{(2)}_m(u) \right|<\rho \kappa\delta \sqrt{\ln(1/\delta) } \right \},  \label{o2}
\eeqn
and denote $\Omega_0(j)=\Omega_1(j) \cap \Omega_2(j)$. 
For simplicity, we suppress the explicit dependence on $j$ in the notations of $\Omega_0$, $\Omega_1$ and $\Omega_2$. Then, the following statements are true.
\begin{lemma}
Let the constant $\rho$ defined in \fr{o2} be such that $0 < \rho< 1/2$. Then on $ \Omega_0$, one has
\be
\frac{1-2\rho}{1-\rho} |g_m(u)| \leq \left|  g^{\delta}_m(u) \right| \leq \frac{1}{1-\rho} |g_m(u)|.  \label{gdelf}
\ee
\end{lemma}
\begin{lemma}
Let $\hat{\beta}_{\omega}$ be defined by \fr{bhat}, and let constant $\kappa$ be such that $\rho \kappa>2$. Then on $\Omega_1$, one has 
\be \label{e2}
\mathbb{E}\left| \hat{\beta}_{\omega}- \beta_{\omega}  \right|^2=O\left( \max\left\{   \eps^2 2^{2j\nu }, \delta^2 2^{2j\nu }   \right\} \right), 
\ee
and on $\Omega_0$, one has
\be \label{e4}
\mathbb{E}\left| \hat{\beta}_{\omega}- \beta_{\omega}  \right|^4=O\left(  \max\left\{   \eps^4 2^{4j\nu +j'}, \delta^4 2^{4j\nu +j'} \right\} \right).
\ee
\end{lemma}
Now, to determine the choice of the thresholds $\lambda^j_{\eps, \delta}$, we define the quantity
\be \label{sg}
S(g^{\delta}_m )_j = \int^1_0\sum\limits_{m \in \Omega_0} \left|  g^{\delta}_m(u) \right|^2du. 
\ee
Then by \fr{wj} and \fr{gm}, we obtain 
\be \label{gdasy}
\int^1_0\sum\limits_{m \in \Omega_0} \left|  g^{\delta}_m(u) \right|^2du \asymp \sum\limits_{m \in \Omega_0} |m|^{-2\nu} \asymp 2^j 2^{-2j\nu }.
\ee
Hence, we choose our thresholds  $\lambda^j_{\eps, \delta}$ of the form
\be \label{lambda}
\lambda_{\epsilon, \delta,j}=2^{j/2} \left[S(g^{\delta}_m )_j \right]^{-1/2} \max\left\{\gamma_1 \epsilon \sqrt{\ln(1/\epsilon)}, \gamma_2 \delta \sqrt{ \ln(1/\delta)} \right\}.
\ee
Finally, choose $J$ and $J'$ such that  
\beqn
J&=&\max\left\{ j \in \mathbb{Z}:  \left[S(g^{\delta}_m )_j \right]^{-1}\leq 2^{-2j}\left(\max\left\{A^{-2}\eps^{2}, A^{-2}\delta^{2} \right\}     \right)^{-1}\right\}, \label{JJ}\\
2^{J'}&\asymp& \left(\max\{ \eps^{2}, \delta^{2} \}  \right)^{-1}.\label{JJ'}
\eeqn
\begin{lemma}
For J  defined in \fr{JJ}, one has as $\epsilon, \delta \rightarrow 0$, 
\be \label{2J}
2^J \asymp \left(\max\{A^{-2}\eps^{2}, A^{-2}\delta^{2}\}\right)^{-\frac{1}{2\nu+1}} .
\ee
\end{lemma}
\begin{remark} Notice that our choices of the threshold $\lambda^j_{\eps, \delta}$, and the finest resolution levels $J$ and $J'$ are completely determined from the data and therefore estimator \fr{fhat} is adaptive. 
\end{remark}
It remains to investigate how the estimator performs both asymptotically and in a finite sample setting. Next, we evaluate the asymptotic minimax lower and upper bounds for the $L^2$-risk. 
\section{Minimax rates of convergence}
To construct minimax lower bounds for the $L^2$-risk, we define the $L^2$-risk over the set  $V$ as 
\be  \label{eq12}
 R(V)=\inf_{\hat{f}} \sup _{f \in V}\EE \left\| \hat{f}-f \right\|^2_2,
\ee
where $\|g\|_2$ is the $L^2$-norm of a function $g$ and the infimum is taken over all possible estimators $\hat{f}$ of $f$. Let us now introduce an assumption that pertains the functional class of $f$.\\
\noindent
{\bf Assumption A.2.} Let $s^*_i=s_i+1/2- 1/p$, $i=1,2$. Assume that $f(t,u)$ belongs to a two-dimensional Besov ball, and its wavelet coefficients $\beta_{\omega}$
satisfy
\be  \label{ball}
B_{p,q}^{s_1,s_2}(A)=\left\{ f \in L^2(U): \left( \sum\limits_{j,j'}2^{( js_1^*+j's_2^*)q}\left( \sum\limits_{k,k'} \mid \beta_{\omega} \mid^p\right)^{\frac{q}{p}} \right)^{\frac{1}{q}} \leq A \right\} .
\ee
The next three statements are true.
\begin{theorem} \label{th:lowerbds} 
Let $\min\{s_1,s_2\}\geq \max\{1/p,1/2 \}$ with $1\leq p,q < \infty$, let $A>0$ and define $s'_i=s_i+1/2-1/p'$, $p'=\min\{p,2\}$  for $i=1,2$. Then under conditions \fr{gm} and \fr{ball}, as $\eps,\delta \rightarrow 0$, simultaneously, 
\be  \label{low}
R_{\eps,\delta}\left(B_{p,q}^{s_1,s_2}(A)\right)\geq CA^2 \left( \frac{ \max\{\eps^2, \delta^2\} }{A^2}  \right)^d,
\ee
where 
\be \label{d}
d=\left\{\begin{array}{ll} 
\frac{2s_2}{2s_2+1}   & \mbox{if}\ \  s_2(2\nu+1)\leq s_1,\\
\frac{2s_1}{2s_1+2\nu+1} & \mbox{if}\ \  (2\nu+1)\left(\frac{1}{p}-\frac{1}{2}\right)< s_1< s_2(2\nu+1), \\
\frac{2s'_1}{2s'_1+2\nu} & \mbox{if}\ \  s_1 \leq(2\nu+1)\left(\frac{1}{p}-\frac{1}{2}\right).
\end{array}\right.
\ee
\end{theorem}
\begin{lemma}
Let $\hat{\beta}_{\omega}$ and $\lambda^j_{\eps, \delta}$ be defined by \fr{bhat} and \fr{lambda} respectively. Define for $\alpha>0$ the set  
\be \label{theta}
\Theta_{\omega,\alpha}=\left\{\theta:  |\hat{\beta}_{\omega}-\beta_{\omega} |> \alpha\lambda^j_{\eps,\delta}  \right\}.
\ee
Then, under $\Omega_0$ and condition \fr{gm}, and as $\varepsilon, \delta \rightarrow 0$, simultaneously, one has 
\be \label{prob}
\Pr(\Theta_{\omega,\alpha})=O \left( \left(\eps^2\right)^{\tau_1} +\left(\delta^2\right)^{\tau_2}  \right),
\ee
where 
\be \label{tau}
\tau_1=\frac{C_1 (1-2\rho)^2 \alpha^2 \gamma^2_1 }{256 \pi^2 C_2} \textit{ and }  \tau_2=\frac{1}{4}\left[ \frac{\alpha\gamma_2(1-\rho)\sqrt{C_1}}{8\pi \sqrt{C_2}  }- \frac{ \sqrt{C_2} 4\pi \rho^2\kappa^2}{\sqrt{C_1}(1-2\rho)}\right]^2 , 
\ee
 and $C_1,  C_2$ appear in \fr{gm}, and $\kappa$, $\rho$ appear in \fr{o2}.
 \end{lemma}
\begin{theorem} \label{th:upperbds}
Let  $\hat{f}(t,u)$ be the wavelet estimator defined in \fr{fhat}, with $J$ and $J'$ given by \fr{JJ} and \fr{JJ'}. Let conditions \fr{gm} and \fr{ball} hold and $\min\{ s_1, s_2 \}\geq \max\{ 1/p, 1/2\}$, with $1\leq p,q< \infty$, and choose $\tau_1$ and $\tau_2$ in \fr{prob}  large enough. Then as $\varepsilon, \delta \rightarrow 0$, simultaneously,
\be \label{upper}
\sup\limits_{f\in B_{p,q}^{s_1,s_2}(A)} \mathbb{E} \left\| \hat{f}- f \right\|^2 \leq CA^2 \left\{ \left( \frac{\eps^2}{A^{2} }  \ln(1/\eps)\right)^d [\ln (1/\eps)]^{d_1} \vee \left(\frac{\delta^2}{A^{2} } \ln^2\delta \right)^d[ \ln(1/\delta) ]^{d_1} \right\},
\ee
where $d$ is defined in \fr{d} and 
\be \label{d1}
d_1=\II\left(s_1 =(2\nu+1)\left(\frac{1}{p}-\frac{1}{2}\right)\right) +\II( s_1=s_2(2\nu+1)).
\ee
\end{theorem}
\begin{remark}	
 Theorems \ref{th:lowerbds} and \ref{th:upperbds} imply that, for the $L^2$-risk, the estimator in \fr{fhat} is asymptotically quasi-optimal within a logarithmic factor of $\varepsilon$ or $\delta$ over a wide range of anisotropic Besov balls $ B^{s_1, s_2}_{p, q}(A)$. 
\end{remark}
\begin{remark}	
The convergence rates are expressed as a maxima between two terms, taking into account both noise sources. Similar behavior has been pointed out in Hoffmann and Reiss~(2008), Vareschi~(2015) and Benhaddou~(2018a, 2018b). These rates depend on a delicate balance between the parameters of the Besov ball, smoothness of the convolution kernel $\nu$ and the noise parameters $\varepsilon$ and $\delta$.
\end{remark}
\begin{remark}	
With $\delta=0$, our convergence rates coincide with those in Benhaddou et al.~(2013), and with $\delta=0$ and $p=2$, our convergence rates match those in Benhaddou~(2017). In addition, if we hold the variable $u$ fixed, then with $\alpha=1$, our rates are comparable to those in Benhaddou~(2018a) in their univariate case, and with $\alpha_{1l}=1$ and $M=1$  our convergence rates match those in Benhaddou~(2018b), in their univariate but multichannel case.\end{remark}


\section{Simulation Study}
In this section, we carry out a limited simulation study in order to investigate the finite sample performance of our estimator. The first step though is to provide the sample equivalent to equations \fr{fconv}, \fr{gdel}, \fr{bhat},  \fr{o1}-\fr{o2} and \fr{lambda}-\fr{J'}. Indeed, apply the Fourier transform to \fr{disc} to obtain
\be \label{fourier}
y_m(u_l)=f_m(u_l)g_m(u_l)+\sigma_1 z^{(1)}_m(u_l),\
g^{\delta}_m(u_l)=g_m(u_l)+\sigma_2 z^{(2)}_m(u_l).
\ee
Then, the discrete version of \fr{bhat}  is given by
\be  \label{sbhat}
\hat{\beta}_{j,k,j',k'}=\frac{1}{M}\sum^M_{l=1}\sum\limits_{m \in W_j}  \frac{y_m(u_l)}{\hat{g}^{\delta}_m(u_l)} \ \eta_{j',k'}(u_l) \  \overline{\psi}_{j,k,m}.
\ee 
In addition, equations \fr{disc} and \fr{conv} are equivalent by setting
\be \label{epsilondelta}
 \eps^2=\frac{\sigma^2_1}{MN}, \ \delta^2=\frac{\sigma^2_2}{MN}.
 \ee
The sets $\Omega_1$ and $\Omega_2$ are now of the form  
\beqn 
\Omega_1&=&\left \{  m\in W_j: \min_{l\leq M}  \left|g^{\delta}_m(u_l)  \right|^2>\frac{\kappa^2 \sigma^2_2}{MN}\ln(MN)  \right \},     \\
\Omega_2&=&\left \{ m\in W_j: \max_{l\leq M}\left|z^{(2)}_m(u_l) \right|^2<\rho^2 \kappa^2 \ln(MN)   \right \}, 
\eeqn
and \fr{sg} has the sample counterpart
\be \label{sg1}
S(g^{\delta}_m )_j =\frac{1}{M} \sum^M_{l=1}\sum\limits_{m \in \Omega_0} \left|g^{\delta}_m(u_l) \right|^2. 
\ee
Finally,  the threshold and the finest resolution levels are given by
\beqn
\lambda_{j,M,N}&=& \gamma 2^{j/2} \left[S(g^{\delta}_m)_j  \right]^{-1/2}  \max\left\{ \left[ \frac{\sigma^2_1 \ln(MN)}{MN} \right]^{1/2},\left[ \frac{\sigma^2_2 \ln^2(MN)}{MN} \right]^{1/2}  \right\}, \\
J&=&\max\left\{ j \in \mathbb{Z}:  \left[S(g^{\delta}_m) _j \right]^{-1}\leq 2^{-2j} \left(  \max\left\{\frac{\sigma^2_1 }{MN} ,\frac{\sigma^2_2 }{MN}\right\}\right)^{-1}  \right\},\label{J}\\
 2^{J'}&=&\left(  \max\left\{\frac{\sigma^2_1 }{MN} ,\frac{\sigma^2_2 }{MN}\right\}\right)^{-1}. \label{J'}
\eeqn
The simulation is implemented based on the above equations. In particular, it is formatted through {\fontfamily{qcr}\selectfont \small MATLAB} using the {\fontfamily{qcr}\selectfont \small Wavelab} toolbox. Similar to Benhaddou et al.~(2013), degree 3 Meyer wavelets family and degree 6 Daubechies wavelets family are utilized for the wavelet transform. We generate our data by equation \fr{fourier} with test kernel $g(t,u)=0.5\exp(-|t|( 1+ ( u-0.5  )^2      ))$, and with various test functions $f(t,u)$ and different combinations of values $N$, $M$ and $\sigma_1$, $\sigma_2$. We use $M=128,\ 256$ and $N=512,\ 1024$. In particular, we generate $f(t,u)$ from $f(t,u)=f_1(t)f_2(u)$ with $f_2(u)$ being a quadratic function $ (u-0.5)^2$ scaled to have a unit norm, and $f_1(t)$ being the routinely used testing functions \textit{Blip}, \textit{Bumps}, \textit{HeaviSine}, and \textit{Doppler} (for details, see Donoho and Johnstone (1994)). We rescale $f(u,t)$ so it has a unit norm.  Graphs of all test functions are presented in Figure 3.1. \\
\noindent As noted in Benhaddou et al.~(2013), our method does not know that $f(t,u)$ was generated by a product of two functions, thus can not take advantage of this prior information. For the choices of $\sigma_1$ and $\sigma_2$, they are determined by the signal-to-noise ratio (SNR),  
\begin{equation*} 
\text{SNR}_1=10\log_{10} \left( \frac{|| f*g ||^2_2}{\sigma^2_1}\right),\ \text{SNR}_2=10\log_{10} \left( \frac{|| g ||^2_2}{\sigma^2_2} \right).
\end{equation*}   
SNR$_1$ is considered for three scenarios,  SNR$_1$=10 dB  (high noise), SNR$_1$= 20 dB (medium noise), and SNR$_1$= 30 dB (low noise). And we consider SNR$_2$=30 dB (low noise). \\
One of the most delicate tasks in implementing the estimation algorithm is the selection
of the finest resolution levels $J$ and $J'$. Theoretically, the choices of $J$ and $J'$ are governed by \fr{J} and \fr{J'}. In a finite sample setting, the algorithm requires them to satisfy 
\begin{equation*}
J\leq \log_2(N)-1, \ J'\leq  \log_2(M)-1. 
\end{equation*}
We choose $J'=\log_2(M)-1$, and empirically investigate the performance of the algorithm by choosing $J$ from the set $\left\{3, 4, 5, \ldots, \log_2(N)-1 \right\}$. The selection of $J$ reported in Table 3.1 (in the brackets) is made when MISE attains the smallest for $J \in \left\{3, 4, 5, \ldots, \log_2(N)-1 \right\}$.\\
\noindent We compute the empirical version of MISE through $N_0=100$ simulation repetitions, where,
\begin{equation*}
\widehat{MISE}(\hat{f},f)= \frac{1}{N_0} \sum_{i=1}^{N_0}\left\| \hat{f}_{i}-f \right\|^2.
\end{equation*}
Table 1 reports the averages of those errors over 100 simulation repetitions together with their standard deviations (in the parentheses).\\
\begin{figure}
\centering
\caption{Testing functions}
\resizebox{6.35in}{4.35in}{\includegraphics{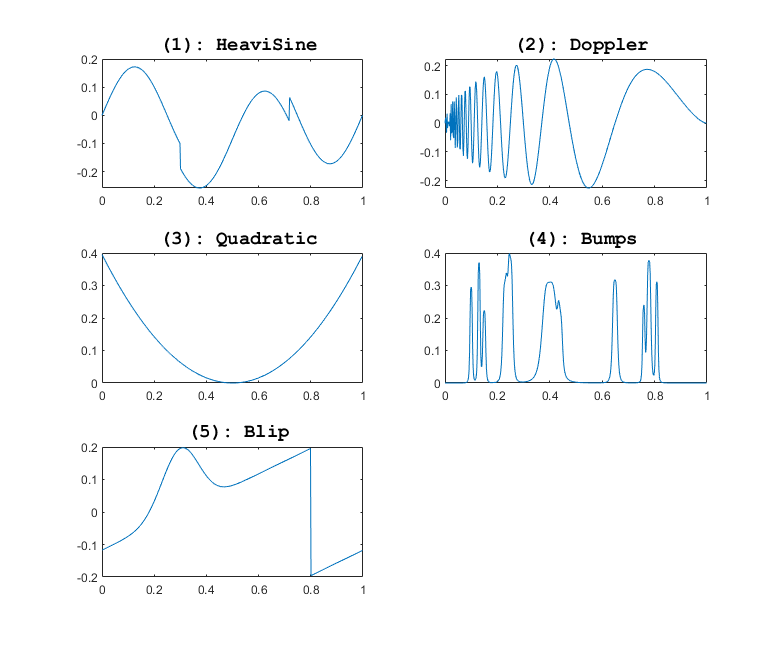} }
\end{figure}
\begin{table}
\centering
\caption{MISE averaged over 100 simulation repetitions}
\label{table}
\resizebox{\columnwidth}{!}{
\begin{tabular}{|c| c| c| c| c| c| c| c| c| c|  c|  c|  c| } 
 \hline
 \multicolumn{1}{|c|}{ N=} &\multicolumn{3}{c|}{512 }&\multicolumn{3}{c|}{1024 }&\multicolumn{3}{c|}{512 }&\multicolumn{3}{c|}{ 1024 }\\
 \hline
  SNR$_1$=& 10dB & 20dB & 30dB & 10dB & 20dB & 30dB & 10dB & 20dB & 30dB & 10dB & 20dB & 30dB \\
\hline
 \multicolumn{1}{|c|}{ } &       \multicolumn{6}{c|}{$f_u=Quadratic, f_t=HeaviSine$  }         &\multicolumn{6}{c|}{$f_u=Quadratic, f_t=Doppler$ }\\
 \hline
M=128 & 0.0085 [3] & 0.0054 [3]   & 0.0051 [3]  & 0.0059 [3]   & 0.0051 [3]   & 0.0027 [4]        & 0.2547 [3]  & 0.1226 [4]   & 0.0480 [5]    & 0.1247 [4]    & 0.0499 [5]     & 0.0480 [5]\\
            & (0.000147) & (0.000013) &(0.000001) & (0.000033)  & (0.000003) & (0.000002)       & (0.000060) & (0.000035)  & (0.000012)   & (0.000078)  & (0.000045)   & (0.000003)\\
 \hline
M=256 & 0.0085 [3] & 0.0053 [4]  & 0.0028 [4]   & 0.0059 [3]   & 0.0033 [4]   & 0.0018 [5]       & 0.1325 [4]   & 0.0562 [5]   & 0.0181 [6]    & 0.0698 [5]    & 0.0316 [6]     & 0.0155 [6]\\
           & (0.000102) & (0.000061) & (0.000005) & (0.000025) & (0.000014)  & (0.000008)      & (0.000260) & (0.000139)  & (0.000073)   & (0.000335)  & (0.000166)    & (0.000018)\\
\hline
\hline
\multicolumn{1}{|c|}{ } &\multicolumn{6}{c|}{$f_u=Quadratic, f_t=Bumps $  }                       &\multicolumn{6}{c|}{$f_u=Quadratic, f_t=Blip$  }\\
\hline
M=128 & 0.3046 [3] & 0.1615 [4]   & 0.0837 [5]    & 0.2966 [3]   & 0.1537 [4]   & 0.0780 [5]     & 0.0609 [3]   & 0.0577 [3]    & 0.0250 [4]   & 0.0586 [3]    & 0.0259 [4]   &  0.0144 [5]\\
            & (0.000498) & (0.000258) & (0.000157)  & (0.000120) & (0.000062) & (0.000038)   & (0.000138) & (0.000016)  & (0.000006)  & (0.000040)  & (0.000021) & (0.000011)\\
\hline 
M=256 & 0.2396 [4] & 0.1470 [5]   & 0.0604 [6]    & 0.1735 [4]   & 0.0939 [5]   & 0.0209 [6]     & 0.0532 [4]   & 0.0277 [4]    & 0.0160 [5]   & 0.0324 [4]    & 0.0199 [5]   & 0.0100 [6]\\
            & (0.001570) & (0.001079) & (0.000614)  & (0.000410) & (0.000260) & (0.000148)   & (0.000563) & (0.000062)  & (0.000035)  & (0.000145)  & (0.000085) &(0.000050)\\
 \hline
\multicolumn{1}{|c|}{ } &\multicolumn{6}{c|}{$f_u=Bumps, f_t=Bumps $  }                                      &\multicolumn{6}{c|}{$f_u=Blip, f_t=Blip$   }\\
\hline
M=128 & 0.3057 [3]   & 0.1625 [4]    & 0.0849 [5]    & 0.2969 [3]     & 0.1539 [4]     & 0.0783 [5]            & 0.0611 [3]    & 0.0577 [3]    & 0.0250 [4]    & 0.0586 [3]   & 0.0259 [4]    & 0.0145 [5]\\
             & (0.000493)  & (0.000268)  & (0.000170)  & (0.000129)   & (0.000076)  & (0.000042)           & (0.000172)  & (0.000017)  & (0.000008)  & (0.000037) & (0.000021)  & (0.000013)\\
\hline	
M=256 & 0.2486 [4]   & 0.1541 [5]    & 0.0662 [6]    & 0.1758 [4]     & 0.0958 [5]    & 0.0224 [6]             & 0.0545 [4]    & 0.0278 [4]    & 0.0161 [5]    & 0.0327 [4]    & 0.0201 [5]    & 0.0102 [6]\\
           & (0.00214)      & (0.001181)  & (0.000616)  & (0.000512)   & (0.000323)  & (0.000152)           & (0.000669)  & (0.000062)  & (0.000041)   & (0.000131) & (0.000088)  & (0.000049)\\	
 \hline
 \hline	
\multicolumn{1}{|c|}{ } &\multicolumn{6}{c|}{$f_u=Blip, f_t=Bumps $   }                              &\multicolumn{6}{c|}{$f_u=Bumps, f_t=Blip$    }\\
\hline	
M=128 & 0.3053 [3]  & 0.1622 [4]    & 0.0844 [5]   & 0.2968 [3]    & 0.1539 [4]   & 0.0782 [5]        & 0.0611 [3]    & 0.0577 [3]   & 0.0251 [4]   & 0.0587 [3]   & 0.0259 [4]    & 0.0145 [5]  \\
           & (0.000441)   & (0.000287)   & (0.000172) & (0.000120)  & (0.000070)  & (0.000037)       & (0.000183)  & (0.000015)  & (0.000008)  & (0.000037)  & (0.000023)  & (0.000013)   \\
\hline
M=256 & 0.2468 [4]  & 0.1526 [5]    & 0.0651 [6]   & 0.1754 [4]   & 0.0953 [5]   & 0.0221 [6]         & 0.0549 [4]    & 0.0278 [4]   & 0.0162 [5]   & 0.0328 [4]   & 0.0202 [5]    & 0.0103 [6]  \\
           & (0.002128)   & (0.001147)  & (0.000695)  & (0.000480)  & (0.000276)  & (0.000162)       & (0.000574)  & (0.000059)  & (0.000041)  & (0.000182)  & (0.000089)  &(0.000061)  \\
\hline
\end{tabular}
}
\end{table}
\noindent The simulation is aimed to study the effect of two components; the sample budget $MN$ and the noise levels $\sigma_i$, $i=1, 2$. The simulation results confirm the theory and are consistent with previous results in the literature. That is, as the sample size increases ($MN$ decreases), the performance of estimation increases. At the same time, as the noise level $R_i$ increases ($\sigma_i$ increases), the performance deteriorates. The simulation results confirm that as the noise level $\sigma_i$ decreases ($R_i$ increases), the convergence rate improves. 
\section{Proofs}
{\bf Proof of Lemma 1.} 
From \fr{conv}, one has that
\be
|g_m(u)| =\left|g^{\delta}_m(u)- \delta z^{(2)}_m(u)\right| 
= \left|g^{\delta}_m(u) \left( 1- \frac{\delta z^{(2)}_m(u)}{g^{\delta}_m(u)}   \right) \right| = \left|  g^{\delta}_m(u) \right| \left| 1- \frac{\delta z^{(2)}_m(u)}{g^{\delta}_m(u)} \right|.
\ee
Therefore,
\be \label{gdsum1}
\left|  g^{\delta}_m(u) \right| = |g_m(u)| \left| \frac{1}{1- \frac{\delta z^{(2)}_m(u)}{g^{\delta}_m(u)}} \right| = |g_m(u)| \left| \sum\limits_{l=0}^{\infty} \left(   \frac{\delta z^{(2)}_m(u)}{g^{\delta}_m(u)} \right)^l    \right|.
\ee
This is because on $\Omega_0$, one has
\be \label{upbnd1}
0<\left|\frac{\delta z^{(2)}_m(u)}{g^{\delta}_m(u)} \right| <\left|\frac{\rho \kappa \delta \sqrt{\ln(1/\delta) }}{\kappa \delta \sqrt{\ln(1/\delta) }} \right|=|\rho| <1     .
\ee
Thus, using \fr{gdsum1} and \fr{upbnd1}, yields
\be
\left|  g^{\delta}_m(u) \right| = |g_m(u)| \left| \sum\limits_{l=0}^{\infty} \left(   \frac{\delta z^{(2)}_m(u)}{g^{\delta}_m(u)} \right)^l    \right|
\leq |g_m(u)|\frac{1}{1-\rho}\label{gdup}.
\ee
On one hand, one has
\be \label{glowb}
|g_m(u)|  > \left|  g^{\delta}_m(u) \right|- \left| \delta z^{(2)}_m(u) \right|
>(1-\rho) \kappa \delta \sqrt{\ln(1/\delta) }.
\ee
On the other hand,
\be \label{gdsum}
\left|  g^{\delta}_m(u) \right|=|g_m(u)+\delta z^{(2)}_m(u)|=\left| g_m(u) \left( 1+ \frac{\delta z^{(2)}_m(u)}{g_m(u)} \right) \right|.
\ee
Combining \fr{o2} and \fr{glowb} with $ 0<\rho<\frac{1}{2}$, one has
\be \label{lbr}
\left| \frac{\delta z^{(2)}_m(u)}{g_m(u)} \right| <\left|\frac{\rho \kappa \delta \sqrt{\ln(1/\delta) }}{(1-\rho) \kappa \delta \sqrt{\ln(1/\delta) }} \right| =\frac{\rho}{1-\rho} <1.
\ee
Now, using \fr{gdsum} and \fr{lbr}, one obtains
\be
|g_m(u)|=\left|  g^{\delta}_m(u) \right| \left| \sum\limits_{l=0}^{\infty} \left(  - \frac{\delta z^{(2)}_m(u)}{g_m(u)} \right)^l\right|
\leq \left|  g^{\delta}_m(u) \right|\frac{1-\rho}{1-2\rho}. \label{gddwn}
\ee
Hence, combining \fr{gdup} and \fr{gddwn} completes the proof. $\Box$\\
\noindent {\bf Proof of Lemma 2.} 
First, an extension of It\^o isometry (see, \O ksendal (2003)) provides that for any function $F(u,t) \in L^2(U)$, one has 
\be 
\EE\left[\int_{0}^{1}\int_{0}^{1} F(u,t)dB(u,t)du   \right]^2=\int_{0}^{1}\int_{0}^{1}F^2(u,t)dtdu.\label{ito}
\ee
Consider the result in \fr{e2}, one has
\be 
\EE \left| \hat{\beta}_{j,k,j',k'}- \beta_{j,k,j',k'} \right|^2 \leq 2(V_1+V_2),
\ee
where
\beqn
V_1 &=&\EE \left|\sum\limits_{m \in W_j} \int_{0}^{1} \left(  \frac{y_m(u)}{\hat{g}^{\delta}_m(u)} -   f_m(u)     \right) \ \eta_{j',k'}(u) du \  \overline{\psi}_{j,k,m} \II(\Omega_1)    \right|^2 ,\\
V_2 &=&\EE \left| \sum\limits_{m \in W_j} \int_{0}^{1} f_m(u)\ \eta_{j',k'}(u) du \  \overline{\psi}_{j,k,m} \II\left(\Omega^c_1\right) \right|^2 .\label{V2}
\eeqn
First, consider $V_2$ in \fr{V2}, using \fr{gm} and Gaussian tail inequality, yields
\begin{align}
\Pr\left(\Omega^c_1\right)=& \Pr\left( \left|g_m(u)+\delta z^{(2)}_m(u)\right|<\kappa \delta \sqrt{\ln(1/\delta) } \right) \leq\Pr\left( |g_m(u)|-|\delta z^{(2)}_m(u)|<\kappa \delta \sqrt{\ln(1/\delta) } \right) \nonumber\\
\leq&\Pr \left(\left| z^{(2)}_m(u)\right|>\frac{1}{\delta} |g_m(u)|-\kappa \sqrt{\ln(1/\delta) } \right) \leq \Pr\left(\left| z^{(2)}_m(u)\right|>\frac{1}{\delta} |g_m(u)|(1-o(1))\right) \nonumber\\
\leq& \Pr\left(\left| z^{(2)}_m(u)\right|>\frac{\sqrt{C_1}|m|^{-\nu}}{\delta} (1-o(1))\right)  
=O\left(   \exp\left\{-\frac{c_1}{2\delta^2 2^{2j\nu}}\right\} \right) =O\left( \delta^2 2^{2j\nu } \right)\label{pc1}.
\end{align}
Thus, using \fr{beta}, \fr{pc1}, and applying Cauchy-Schwarz inequality, one has
\beqn
V_2&=&\EE \left|  \int_{0}^{1}  \int_{0}^{1} f(u,t)\psi_{j,k}(t) \eta_{j',k'}(u) \ dt du\  \II\left(\Omega^c_1\right) \right|^2 \nonumber\\
&\leq&\int_{0}^{1}   \int_{0}^{1}  \left| \psi_{j,k}(t) \eta_{j',k'}(u)   \right|^2 dtdu\int_{0}^{1}   \int_{0}^{1}  \left| f(u,t)   \right|^2 dtdu \Pr \left(\Omega^c_1\right)= O\left( \delta^2 2^{2j\nu} \right).  \label{v2}
\eeqn
And $V_1$ can be further partitioned as
\begin{align}
V_1=&\EE \left|\sum\limits_{m \in W_j} \int_{0}^{1} \left(  \frac{\eps z^{(1)}_m(u)-\delta f_m(u)z^{(2)}_m(u)}{g^{\delta}_m(u)}  \right) \ \eta_{j',k'}(u) du \  \overline{\psi}_{j,k,m} \II(\Omega_1)  \right|^2  
\leq 2\left(\EE_1+\EE_2 \right),
\end{align} 
where
\beqn
\EE_1&=&\EE \left|\sum\limits_{m \in W_j} \int_{0}^{1}\left( \frac{\eps z^{(1)}_m(u)-\delta f_m(u)z^{(2)}_m(u)}{g^{\delta}_m(u)}  \right) \ \eta_{j',k'}(u) du \  \overline{\psi}_{j,k,m} \II(\Omega_1) \II(\Omega_2) \right|^2, \\
\EE_2&=&\EE \left|\sum\limits_{m \in W_j} \int_{0}^{1}\left( \frac{\eps z^{(1)}_m(u)-\delta f_m(u)z^{(2)}_m(u)}{g^{\delta}_m(u)}  \right) \ \eta_{j',k'}(u) du \  \overline{\psi}_{j,k,m} \II(\Omega_1) \II\left(\Omega^c_2\right) \right|^2.
\eeqn
And $\EE_1$ and $\EE_2$ can be further partitioned as $\EE_1\leq 2( \EE_{11}+\EE_{12})$ and $\EE_2\leq 2( \EE_{21}+\EE_{22})$, where
\beqn
\EE_{11}&=&\EE \left|\sum\limits_{m \in W_j} \int_{0}^{1}\left( \frac{\eps z^{(1)}_m(u)}{g^{\delta}_m(u)}  \right) \ \eta_{j',k'}(u) du \  \overline{\psi}_{j,k,m}   \II(\Omega_1) \II(\Omega_2) \right|^2, \\
\EE_{12}&=&\EE \left|\sum\limits_{m \in W_j} \int_{0}^{1}\left( \frac{\delta f_m(u)z^{(2)}_m(u)}{g^{\delta}_m(u)}  \right) \ \eta_{j',k'}(u) du \  \overline{\psi}_{j,k,m}  \II(\Omega_1) \II(\Omega_2) \right|^2 ,\\
\EE_{21} &=& \EE \left|\sum\limits_{m \in W_j} \int_{0}^{1}\left( \frac{\eps z^{(1)}_m(u)}{g^{\delta}_m(u)}  \right) \ \eta_{j',k'}(u) du \  \overline{\psi}_{j,k,m}  \II(\Omega_1) \II\left(\Omega^c_2\right) \right|^2,\label{E_21} \\
\EE_{22}& =&\EE \left|\sum\limits_{m \in W_j} \int_{0}^{1}\left( \frac{\delta f_m(u)z^{(2)}_m(u)}{g^{\delta}_m(u)}  \right) \ \eta_{j',k'}(u) du \  \overline{\psi}_{j,k,m} \II(\Omega_1) \II\left(\Omega^c_2\right) \right|^2. 
\eeqn
Consider $\EE_{11}$, take into account of \fr{wj}, \fr{gdelf}, \fr{ito} and the fact that $|\psi_{j,k,m}|\leq 2^{-j/2}$, one has 
\beqn
\EE_{11}
&=&\eps^2\EE \left|\sum\limits_{m \in W_j} \int_{0}^{1}\left( \frac{  \int_{0}^{1} e_m(t)z^{(1)}(u,t) dt}{g^{\delta}_m(u)}  \right) \ \eta_{j',k'}(u) du \  \overline{\psi}_{j,k,m} \II(\Omega_1) \II(\Omega_2) \right|^2 \nonumber\\
&=&\eps^2\EE \left| \int_{0}^{1} \int_{0}^{1}  \sum\limits_{m \in W_j} \frac{   e_m(t)}{g^{\delta}_m(u)}   \ \eta_{j',k'}(u)\  \overline{\psi}_{j,k,m}    dB^{(1)}(u,t) du  \right|^2 \nonumber\\
&=&\eps^2\int_{0}^{1} \int_{0}^{1} \left[ \sum\limits_{m \in W_j} \frac{  e_m(t)}{g^{\delta}_m(u)}   \ \eta_{j',k'}(u)\  \overline{\psi}_{j,k,m} \right]^2   dt du  \nonumber\\
&=&\eps^2\int_{0}^{1} \int_{0}^{1}\sum\limits_{m_1  }\sum\limits_{m_2  }  \frac{  e_{m_1}(t)e_{m_2}(t) }{g^{\delta}_{m_1}(u)\overline{g}^{\delta}_{m_2}(u)} \overline{\psi}_{j,k,m_1}\psi_{j,k,m_2} \left| \eta_{j',k'}(u) \right|^2 dt du\nonumber\\ 
&=& \eps^2 \int_{0}^{1} \int_{0}^{1}  \sum\limits_{m \in W_j} \frac{ | \eta_{j',k'}(u)|^2 | \overline{\psi}_{j,k,m} |^2}{|g^\delta_m(u)|^2}   \ dt du =O\left( \eps^2 2^{2j\nu}\right) . \label{E11}
\eeqn
Note that in the double summation above, all terms involving $m_1 \neq m_2$ vanish due to the fact that  $ \int_{0}^{1} e_{m_1}(t)e_{m_2}(t) dt=0$.
\noindent Using the same arguments, one has that 
\begin{equation}
\EE_{12} = O\left( \delta^2 2^{2j\nu}\right). \label{eq: 2.26}
\end{equation}
To evaluate $\EE_{21}$ in \fr{E_21},  notice that applying Gaussian tail inequality on $\Omega^c_2$, one has
\be
\PP\left(\Omega^c_2\right)= \PP \left(  |\delta z^{(2)}_m(u)|^2\geq \rho^2 \kappa^2 \delta^2 \ln(1/\delta)  \right) \leq \exp\left\{-\frac{1}{2}\rho^2\kappa^2\ln(1/\delta)\right\}=\delta^{\frac{1}{2}\rho^2\kappa^2}. \label{O2c}
\ee
Thus, using Cauchy-Schwarz inequality and taking into account of the condition that $\rho^2 \kappa^2>8$ and using \fr{gdelf} and \fr{O2c}, one has
\begin{align}
\EE_{21}\leq&\EE\int_{0}^{1} \left|\sum\limits_{m \in W_j} \left( \frac{\eps z^{(1)}_m(u)}{g^{\delta}_m(u)}  \right) \ \eta_{j',k'}(u)   \overline{\psi}_{j,k,m} \II(\Omega_1) \II\left(\Omega^c_2\right) \right|^2 du \nonumber\\
=&\eps^2\EE \left[\int_{0}^{1} \sum\limits_{m=m' \in W_j }  \frac{ z^{(1)}_m(u_1)}{g^{\delta}_m(u_1)}\II(\Omega_1)\frac{ z^{(1)}_{m'}(u_2)}{\overline{g^{\delta}}_{m'}(u_2)}\II(\Omega'_1)\eta_{j',k'}(u_1)\eta_{j',k'}(u_2) \overline{\psi}_{j,k,m} {\psi}_{j,k,m'} \II\left(\Omega^c_2\right)\II(\Omega'^c_2)du \right] \nonumber\\
=&\eps^2  \sum\limits_{m=m' \in W_j} \left[\int_{0}^{1}  |\eta_{j',k'}(u)|^2 | \overline{\psi}_{j,k,m}|^2du\ \EE\left( \frac{ z^{(1)}_m \II(\Omega_1)}{g^{\delta}_m} \frac{ z^{(1)}_{m'}\II(\Omega'_1)}{\overline{g^{\delta}}_{m'}} \II\left(\Omega^c_2\right)\II(\Omega'^c_2)   \right) \right] \nonumber\\
\leq&\eps^2  \sum\limits_{m=m' \in W_j} \int_{0}^{1}  |\eta_{j',k'}(u)|^2 | \overline{\psi}_{j,k,m}|^2du\ \sqrt{\EE\left( \frac{ z^{(1)}_m \II(\Omega_1)}{g^{\delta}_m} \frac{ z^{(1)}_{m'} \II(\Omega'_1)}{\overline{g^{\delta}}_{m'}} \right)^2 \EE\left(\II\left(\Omega^c_2\right)\II(\Omega'^c_2) \right)^2 } \nonumber\\
=&\eps^2 \sum\limits_{m=m' \in W_j} \int_{0}^{1}  |\eta_{j',k'}(u)|^2 | \overline{\psi}_{j,k,m}|^2du\ \sqrt{\EE\left( \frac{ z^{(1)}_m\II(\Omega_1)}{g^{\delta}_m}  \right)^2 \EE\left(\frac{ z^{(1)}_{m'} \II(\Omega'_1) }{\overline{g^{\delta}}_{m'}} \right)^2 \PP\left(\Omega^c_2\right)\PP(\Omega'^c_2) } \nonumber\\
\leq&C\eps^2 \sum\limits_{m=m' \in W_j} \int_{0}^{1}  |\eta_{j',k'}(u)|^2 | \overline{\psi}_{j,k,m}|^2du \frac{1}{ \kappa^2 \delta^2 \ln(1/\delta)} \delta^{\frac{1}{2}\rho^2\kappa^2} =O\left( \max\left\{\eps^2 ,\delta^2  \right\}\right) \label{E21}.
\end{align}
Similarly, one has
\be 
\EE_{22}=O\left( \max\left\{\eps^2 ,\delta^2  \right\}\right).\label{E22}
\ee
Combining results \fr{E11} - \fr{E22},  yields \fr{e2}.\\
\noindent Next, we prove the result in \fr{e4}. On $\Omega_0$, we have 
\be 
\EE \left| \hat{\beta}_{j,k,j',k'} - \beta_{j,k,j',k'}  \right|^4  \leq 2^3 \left(T_1+T_2+T_3\right), \label{beta4}
\ee
where
\beqn 
T_1& =& \EE \left|\sum\limits_{m \in W_j} \int_{0}^{1}   \frac{\eps z^{(1)}_m(u) }{g^{\delta}_m(u)}  \eta_{j',k'}(u) du \  \overline{\psi}_{j,k,m}\II(\Omega_0) \right|^4 ,\\
T_2&=&\EE \left|\sum\limits_{m \in W_j} \int_{0}^{1} \frac{\delta f_m(u)z^{(2)}_m(u)}{g^{\delta}_m(u)}   \eta_{j',k'}(u) du \  \overline{\psi}_{j,k,m}\II(\Omega_0) \right|^4.\\
T_3&=&\EE\left| \sum\limits_{m \in W_j} \int_{0}^{1} f_m(u)\ \eta_{j',k'}(u) du \ \overline{\psi}_{j,k,m}   \II\left(\Omega^c_0\right) \right|^4 .
\eeqn
First, consider $T_3$,  one has
\be
 \PP\left(\Omega^c_0\right) \leq  \PP\left(\Omega^c_1\right)+ \PP\left(\Omega^c_2\right). \label{o0}
\ee
Thus,  under the condition that $\rho^2 \kappa^2 > 8$, using \fr{beta} and applying Cauchy-Schwarz inequality and \fr{o0}, one has
\beqn
T_3&=&\EE\left|  \int_{0}^{1}  \int_{0}^{1} f(u,t)\psi_{j,k}(t) \eta_{j',k'}(u) \ dt du\   \II\left(\Omega^c_0\right)  \right|^4  \nonumber\\
&\leq&\EE\left[ \left( \int_{0}^{1} \int_{0}^{1} \left| f(u,t) \right|^2 dt du\  \II\left(\Omega^c_0\right)  \right)^2   \left( \int_{0}^{1} \int_{0}^{1} \left| \psi_{j,k}(t) \eta_{j',k'}(u)\right|^2 dtdu     \right)^2 \right] \nonumber\\ 
&=& \PP\left(\Omega^c_0\right)=O\left(   \delta^4 2^{4j\nu}  \right).
\eeqn
Using  \fr{zm}, \fr{wj}, \fr{gm},  \fr{gdelf} and applying Cauchy-Schwarz inequality and Minkowski's inequality and note that $|\psi_{j,k,m}|\leq 2^{-j/2}$, one has
\beqn
T_{1}&=&\epsilon^4  \EE  \left| \sum\limits_{m \in W_j} \int_{0}^{1} \frac{\int_{0}^{1} e_m(t)z^{(1)}(u,t) dt}{g^{\delta}_m(u)} \eta_{j',k'}(u) du \  \overline{\psi}_{j,k,m}\II(\Omega_0) \right|^4 \nonumber\\
&=&\epsilon^4  \EE \left| \int_{0}^{1}\int_{0}^{1}  \sum\limits_{m \in W_j} \frac{e_m(t) \overline{\psi}_{j,k,m} \II(\Omega_0)}{g^{\delta}_m(u)}  \eta_{j',k'}(u) dt du  \right|^4 \nonumber\\
 &\leq&\epsilon^4  \EE\left[ \int_{0}^{1}\int_{0}^{1}  \left|  \sum\limits_{m \in W_j}   \frac{e_m(t)\overline{\psi}_{j,k,m}\II(\Omega_0) }{g^{\delta}_m(u)}  \right|^2  dtdu \int_{0}^{1}\int_{0}^{1}\left| z^{(1)}(u,t) \eta_{j',k'}(u)\right|^2dtdu   \right]^2  \nonumber\\
 &=&\epsilon^4  \EE\left[ \int_{0}^{1}\int_{0}^{1}  \sum\limits_{m \in W_j}   \frac{  \left| \overline{\psi}_{j,k,m} \right|^2\II(\Omega_0) }{ \left| g^{\delta}_m(u)\right|^2}    dtdu \int_{0}^{1}\int_{0}^{1}\left| z^{(1)}(u,t) \eta_{j',k'}(u)\right|^2dtdu   \right]^2  \nonumber\\
 & \leq &\epsilon^4 2^{4j\nu}   \int_{0}^{1}\int_{0}^{1}\left| \eta_{j',k'}(u)\right|^4 \EE\left(\left| z^{(1)}(u,t)\right|^4\II(\Omega_0) \right)  dtdu  \nonumber\\
   & = &O\left(  \epsilon^4 2^{4j\nu}2^{j'} \right).
\eeqn
Applying similar arguments to $T_{2}$, one has 
\be 
T_{2}=O\left( \delta^4  2^{4j\nu}  2^{j'}    \right).\label{T2}
\ee
Combine the results \fr{beta4} to \fr{T2} completes the proof. $\Box$\\
{\bf Proof of Theorem 1.} 
In order to prove the theorem, we investigate the cases $\varepsilon=0$ and $\delta=0$, separately. In each case, we consider the lower bounds obtained when the worst functions $f$ (i.e. the hardest to estimate) are uniformly spread over the unit interval; the dense-dense case, and when the worst functions $f$ are sparse; the sparse-dense case. Lemma $A.1$ of Bunea et al.~(2007) is then applied to find such lower bounds using conditions \fr{gm} and \fr{ball}. To complete the proof, one chooses the highest of the lower bounds by comparing the outcomes based on the two cases $\varepsilon=0$ and $\delta=0$. Let $G$ be the class of functions $g(\cdot,\cdot)$ such that 
\begin{equation} 
G=\left\{g(t,x):  C_1 |m|^{-2\nu} \leq |g_m(x)|^2 \leq C_2 |m|^{-2\nu},\ t,x\in [0,1]  \right\},
\end{equation}
and define functions $g_{j, j'}$ as
\be 
g_{j,j'}(t,x)=c_0 2^{-j(\nu-1/2)}\psi_{j,k}(t)\eta_{j',k'}(x).
\ee
Notice that $g_{j,j'}(t,x) \in G$.\\
First, we consider the case when $\eps=0$. Then \fr{conv} becomes 
\be  \label{}
y(t,u)=\int_{0}^{1}f(x,u)g(t-x,u)dx,\
g^{\delta}(t,u)=g(t,u)+\delta z^{(2)}(t,u).
\ee
Define
\be
h^{\delta}(t,u)=\int_{0}^{1}f(x,u)g(t-x,u)dx+ \delta z(t,u),
\ee
where $z(t,u)$ is a white Gaussian noise.\\
\underline{The dense-dense case}. Let $\omega$ be the matrix with elements $\omega_{k,k'}=\{0, 1\}$, $k=0,\ldots, 2^j-1$, $k'=0,\ldots, 2^{j'}-1$. Denote the set of all possible $\omega$ by $\Omega$ and note that $\omega$ have $N=2^{j+j'}$ elements, so $card(\Omega)=2^N$. Define the functions 
\be \label{78}
f_{j,j'}(t,u)=\gamma_{j,j'}\sum\limits_{k=0}^{2^j-1}\sum\limits_{k'=0}^{2^{j'}-1}\omega_{k,k'}\psi_{j,k}(t)\eta_{j',k'}(u),
\ee
 and note that with the choice $\gamma_{j,j'} = A2^{-j\left(s_1+1/2\right)-j'\left(s_2+ 1/2 \right)}$,  $f_{j,j'} \in B_{p,q}^{s_1,s_2}(A)$. If $\hat{f}_{jj'}$ is of the from \fr{78}, with $\hat{\omega}_{kk'}\in \Omega$ instead of $\omega_{kk'}$, then the $L^{2}$-norm of the difference is 
\be \label{80}
\left\| \hat{f}_{j,j'}(t,u)-f_{j,j'}(t,u)\right\|^2=\gamma^2_{j,j'}\sum\limits_{k=0}^{2^j-1}\sum\limits_{k'=0}^{2^{j'}-1}\II(\hat{\omega}_{k,k'}\neq \omega_{k,k'})=\gamma^{2}_{j,j'}H\left( \hat{\omega}, \omega \right),
\ee
where $H\left( \hat{\omega}, \omega \right)= \sum\limits_{k=0}^{2^j-1}\sum\limits_{k'=0}^{2^{j'}-1}\II(\hat{\omega}_{k,k'}\neq \omega_{k,k'})$ is the Hamming distance between the binary sequences $\hat{\omega}$ and $\omega$. The Varshamov-Gilbert lower bound (see Tsybakov (2009), page 104) indicates that one can choose a subset $\Omega' \subset \Omega$ of cardinality $card(\Omega')\geq 2^{N/8}$ and $H(\hat{\omega},\omega)\geq N/8$ for any $\hat{\omega},\omega \in \Omega'$.
Consequently, one has
\be
\left\| \hat{f}_{j,j'}(t,u)-f_{j,j'}(t,u)\right\|^2 \geq \gamma^2_{j,j'}\frac{2^{j+j'}}{8}.
\ee
Let $W(t,u)$ and $\hat{W}(t,u)$ be two Wiener sheets on $U$. Let $\hat{z}(t,u)=\frac{1}{\delta}(\hat{f}_{j,j'}* g_{j,j'} )(t,u)+z(t,u)$, where $z(t,u)=\overset{.}{W}(t,u)$ and $\hat{z}(t,u)=\overset{.}{\hat{W}}(t,u)$. (i.e., $W(t,u)$ and $\hat{W}(t,u)$ are the primitives of  $z(t,u)$ and $\hat{z}(t,u)$, respectively.)
Then, assuming that $  \frac{1}{\delta^2}\int_{0}^{1}\int_{0}^{1} \left[(\hat{f}_{j,j'}-f_{j,j'})* g_{j,j'} \right]^2dtdu < \infty  $, by the multiparameter Girsanov formula (see, e.g., Dozzi (1989), page 89), we get
\be
\mathbb{E}\left[ \ln\left( \PP_f/\PP_{\hat{f}} \right)\right]_{\PP_f} =\mathbb{E}\left[ \frac{1}{2\delta^2}\int_{0}^{1}\int_{0}^{1} \left[(\hat{f}_{j,j'}-f_{j,j'})* g_{j,j'} \right]^2dtdu   \right]. 
\ee
Applying Plancherel's Identity, and taking into account that $|\psi_{j,k,m}|\leq 2^{-j/2}$ on $W_j$ and $|\hat{\omega}_{k,k'}- \omega_{k,k'}| \leq1$, the Kullback divergence can be written as
\be
K(\PP_f,\PP_{\hat{f}})\leq \frac{1}{2\delta^2} \gamma^2_{j,j'}\sum\limits_{k=0}^{2^j-1}\sum\limits_{k'=0}^{2^{j'}-1}\sum\limits_{m\in W_j}\int_{0}^{1}\psi^2_{j,k,m}\eta^2_{j',k'}(u) g^2_{j,j',m}(u)du\leq \frac{C}{2\delta^2}\gamma^2_{j,j'}2^{j+j'}2^{-2j\nu }.
\ee
Now, in order to apply Lemma A.1 of Bunea et al.~(2007) we must choose  $j$ and $j'$ such that
\be \label{s1}
\frac{1}{2\delta^2}A^22^{-2js_1-2j's_2}2^{-2j\nu }\leq C2^{j+j'}/16.
\ee
Therefore, we need to find a combination $\{j, j'\}$ which is the solution to the following optimization problem
\be
2js_1 + 2j's_2 \xrightarrow{\min} j(2s_1+2\nu +1)+j'(2s_2 + 1) \geq \log_2(CA^2/\delta^2), \ j, j \geq 0.
\ee
 It can be shown that the solution is $(j,j')=\left((2s_1+2\nu+1)^{-1}\log_2(CA^2/\delta^2),0\right),$ if $s_2(2\nu+1)>s_1$, and $(j,j')=\left(0,(2s_2+1)^{-1}\log_2(CA^2/\delta^2)\right),$ if $s_2(2\nu+1)\leq s_1$. Consequently,  applying Lemma A.1 of Bunea et al.~(2007), the corresponding lower bounds are 
\beqn  \label{s11}
s^2=
\begin{cases}
CA^2\left(\frac{\delta^2}{A^2} \right)^{\frac{2s_1}{2s_1+2\nu+1}} \text{ if } s_2(2\nu+1)>s_1, \\
CA^2\left(\frac{\delta^2}{A^2} \right)^{\frac{2s_2}{2s_2+1}} \text{\ \ \ \  \ if } s_2(2\nu+1)\leq s_1.
\end{cases}
\eeqn
\underline{The sparse-dense case}. Let $\omega$ be the vector with elements $\omega_{k'}=\{0, 1\}$, $k'=0,\ldots, 2^{j'}-1$. Let $\Omega$ be the set of all possible $\omega$, so that $\omega$ have $N=2^{j'}$ elements, and $card(\Omega)=2^N$. Define the functions  
\be \label{85}
f_{j,j'}(t,x)=\gamma_{j,j'}\sum\limits_{k'=0}^{2^{j'}-1}\omega_{k'}\psi_{j,k}(t)\eta_{j',k'}(x),
\ee
and note that with the choice $\gamma_{j,j'} = A2^{-js^*_1-j'\left(s_2+ 1/2 \right)}$, $f_{j,j'} \in B_{p,q}^{s_1,s_2}(A)$. Following similar reasoning, it is easy  to check that by  Lemma A.1 of Bunea et al.~(2007), we must choose $j,j'$ such that
\be  \label{s2}
\frac{1}{2\delta^2}A^22^{-2js^*_1-2j's_2}2^{-2j\nu }\leq C2^{j'}/16.
\ee
Therefore, we need to find a combination $\{j, j'\}$ which is the solution to the following optimization problem
\be
2js^*_1 + 2j's_2 \xrightarrow{\min} j(2s^*_1+2\nu)+j'(2s_2+1)\geq \log_2(CA^2/\delta^2), \ j, j' \geq 0.
\ee
It is easy to check that the solution is $(j,j')=\left((2s^*_1+2\nu)^{-1}\log_2(CA^2/\delta^2),0\right),$ if $2s_2\nu>s^*_1$, and $(j,j')=\left(0,(2s_2+1)^{-1}\log_2(CA^2/\delta^2)\right),$ if $2s_2\nu\leq s^*_1$. Thus, the lower bounds are
\beqn \label{s22}
s^2=
\begin{cases}
CA^2\left(\frac{\delta^2}{A^2} \right)^{\frac{2s^*_1}{2s^*_1+2\nu}} \text{ if } 2s_2\nu>s^*_1, \\
CA^2\left(\frac{\delta^2}{A^2} \right)^{\frac{2s_2}{2s_2+1}} \text{\ \ if } 2s_2\nu\leq s^*_1.
\end{cases}
\eeqn
Now, consider case when $\delta=0$. The proof will be the same as that of the lower bounds obtained in Benhaddou et al.~(2013) and therefore we skip it.  To complete the proof, choose the highest of the lower bounds by comparing the outcomes based on the two cases $\varepsilon=0$ and $\delta=0$. $\Box$
{\bf Proof of Lemma 4.} 
Now we consider the probability 
\beqn
\Pr(\Theta_{\omega,\alpha}) &=& \Pr \left(\left| \sum\limits_{m \in W_j} \int_{0}^{1}  \frac{\eps z^{(1)}_m(u)}{\hat{g}^{\delta}_m(u)} \ \eta_{j',k'}(u) \II(\Omega_1)\II(\Omega_2)du \  \overline{\psi}_{j,k,m}  \right| > \frac{\alpha}{2}\lambda^j_{\eps,\delta} \right) \nonumber\\
&&\>+\Pr \left(\left| \sum\limits_{m \in W_j} \int_{0}^{1} \frac{\delta f_m(u)z^{(2)}_m(u)}{\hat{g}^{\delta}_m(u)} \ \eta_{j',k'}(u) du \II(\Omega_1)\II(\Omega_2)\  \overline{\psi}_{j,k,m}  \right| > \frac{\alpha}{2}\lambda^j_{\eps,\delta} \right)  \nonumber\\
&:=& {\Pr}_1+{\Pr}_2 . 
\eeqn
To evaluate ${\Pr}_1$, observe that, since $z^{(1)}$ and $z^{(2)}$ are independent, the conditional distribution of $\aleph_1=\sum\limits_{m \in \Omega_0} \int_{0}^{1}  \frac{\eps z^{(1)}_m(u)}{\hat{g}^{\delta}_m(u)} \ \eta_{j',k'}(u)du \  \overline{\psi}_{j,k,m}$ given  $z^{(2)}_m$ is a zero-mean Gaussian with variance
\be
S^2_1\leq  \frac{4 \pi}{C_1}\left( \frac{1-\rho}{1-2\rho} \right)^2  \left(\frac{8\pi}{3}\right)^{2\nu} 2^{2j\nu } \eps^2.
\ee
Thus, applying Gaussian tail inequality, one has
\be
{\Pr}_1 =O \left( \left(\eps^2\right)^{\frac{C_1 (1-2\rho)^2 \alpha^2 \gamma^2_1 }{256 \pi^2C_2 }}\right).
\ee
It can be easily verified that 
\be 
\frac{\delta f_m(u)z^{(2)}_m(u)}{g_m(u)}=\frac{\delta f_m(u)z^{(2)}_m(u)}{g^{\delta}_m(u)}+\frac{\delta^2 f_m(u)(z^{(2)}_m(u))^2}{g^{\delta}_m(u)g_m(u)}.
\ee
Therefore, for ${\Pr}_2$, we have
\beqn
{\Pr}_2&=&\Pr \left(\left| \sum\limits_{m \in \Omega_0} \int_{0}^{1} \left( \frac{\delta f_m(u)z^{(2)}_m(u)}{g_m(u)} -\frac{\delta^2 f_m(u)(z^{(2)}_m(u))^2}{g^{\delta}_m(u)g_m(u)} \right) \eta_{j',k'}(u) du \  \overline{\psi}_{j,k,m}  \right| > \frac{\alpha}{2}\lambda^j_{\eps,\delta} \right)  \nonumber\\
&\leq&\Pr  \left( 
\left| \sum\limits_{m \in \Omega_0} \int_{0}^{1} \frac{\delta f_m(u)z^{(2)}_m(u)}{g_m(u)}  \eta_{j',k'}(u) du \  \overline{\psi}_{j,k,m}  \right|  \right. \nonumber\\
&& \qquad \left. + \left| \sum\limits_{m \in \Omega_0} \int_{0}^{1} \frac{\delta^2 f_m(u)(z^{(2)}_m(u))^2}{g^{\delta}_m(u)g_m(u)} \eta_{j',k'}(u) du \  \overline{\psi}_{j,k,m}  \right| > \frac{\alpha}{2}\lambda^j_{\eps,\delta} 
\right). \label{p2}
\eeqn
Denote the second term in \fr{p2} and within absolute values by $\aleph_3$. Then, for $|f_m(u)|<1$, one has
\beqn \label{second}
\left| \aleph_3  \right|& \leq& \sum\limits_{m \in \Omega_0} \int_{0}^{1} \frac{ |f_m(u)|\delta^2|(z^{(2)}_m(u))^2|}{\left|  g^{\delta}_m(u) \right||g_m(u)|} | \eta_{j',k'}(u)| du \  |\overline{\psi}_{j,k,m}|  \nonumber\\
&\leq & \frac{1-\rho}{1-2\rho} \frac{4\pi  \rho^2\kappa^2\delta^2 \ln \left(1/\delta \right) }{C_1}  \left(\frac{8\pi}{3}\right)^{2\nu}2^{2j\nu } 2^{\frac{j}{2}}.
\eeqn
Now, consider the first term in \fr{p2} and within absolute values and denote it by $\aleph_2$, and notice that this is a zero-mean Gaussian random variable with variance
\be
S^2_2\leq \frac{4\pi \delta^2}{C_1}  \left(\frac{8\pi}{3}\right)^{2\nu}2^{2j\nu } .
\ee
Consequently, by the Gaussian Tail inequality, ${\Pr}_2$ yields
\be \label{p2bound}
{\Pr}_2 =O \left( \exp\left\{ -\frac{1}{2} \left[ \frac{\alpha\gamma_2(1-\rho)\sqrt{C_1}}{8\pi \sqrt{C_2}  }- \frac{ \sqrt{C_2} 4\pi \rho^2\kappa^2}{\sqrt{C_1}(1-2\rho)}\right] ^2  \ln \left(1/\delta \right)    \right\} \right) =O \left(   \left(  \delta^2   \right)^{ \tau_2  }  \right).
\ee $\Box$\\
{\bf Proof of Theorem 2.}
Denote
\be \label{chijj}
\chi_{\eps,\delta, A}=\max \left\{A^{-2} \eps^2{\ln (1/\eps)}, A^{-2} \delta^2\ln^2( 1/\delta) \right\}, \ \ 2^{j_0}=(\chi_{\eps, \delta, A})^{-\frac{d}{2s'_1}}, \ \ 2^{j'_0}=(\chi_{\eps, \delta, A})^{-\frac{d}{2 s'_2}}. 
\ee
and observe that with $J$ and $J'$ given by \fr{JJ} and \fr{JJ'}, the estimation error can be decomposed into the sum of four components as follows
\be
\mathbb{E} \| \hat{f}(u,t)- f(u,t) \|^2 \leq \mathbb{E}_1 +\mathbb{E}_2 +\mathbb{E}_3+\mathbb{E}_4,  \label{e1234}
\ee
where 
\beqn
\mathbb{E}_1&=&\sum\limits_{k=0}^{2^{m_0}-1} \sum\limits_{k'=0}^{2^{m'_0}-1}\mathrm{Var}( \hat{\beta}_{m_0,k,m'_0,k'}),  \\
\mathbb{E}_2&=&\sum\limits_{\omega \in \Omega(J,J')} \mathbb{E}\left[|\hat{\beta}_{\omega}- \beta_{\omega} |^2\II\left(|\hat{\beta}_{\omega} |>\lambda^j_{\eps,\delta}\right) \right], \\
\mathbb{E}_3&=&\sum\limits_{\omega \in \Omega(J,J')} \mathbb{E}\left[| \beta_{\omega}|^2\II\left(|\hat{\beta}_{\omega} |<\lambda^j_{\eps,\delta}\right) \right],\\
\mathbb{E}_4&=&\left(\sum\limits_{j=m_0-1}^{J-1} \sum\limits_{j'=J'}^\infty+\sum\limits_{j=J}^\infty \sum\limits_{j'=m_0'-1}^{J'-1} +\sum\limits_{j=J}^\infty \sum\limits_{j'=J'}^\infty \right) \sum\limits_{k=0}^{2^j-1} \sum\limits_{k'=0}^{2^{j'}-1}  \beta^2_{\omega} .
\eeqn
Combining $\mathbb{E}_1$ and $\mathbb{E}_4$, using \fr{e2} in $\mathbb{E}_1$ and \fr{ball} with \fr{J'} and \fr{2J} in $\mathbb{E}_4$, one has 
\be \label{ee14}
\mathbb{E}_1+\mathbb{E}_4=O\left(\max\{\eps^2, \delta^2 \}+ \left( \sum\limits_{j=J}^\infty + \sum\limits_{j'=J'}^\infty \right) A^2 2^{-2(js'_1+j's'_2)} \right)=O\left( A^2\chi^d_{\eps,\delta, A}\right).
\ee
Notice that 
\beqn
\II\left(|\hat{\beta}_{\omega} |<\lambda^j_{\eps,\delta} \right)&\leq& \II\left(|\hat{\beta}_{\omega} -\beta_{\omega} |> \frac{1}{2}\lambda^j_{\eps,\delta}\right) + \II\left(|\beta_{\omega} |<\frac{3}{2}\lambda^j_{\eps,\delta}\right),\\
\II\left(|\hat{\beta}_{\omega} |>\lambda^j_{\eps,\delta} \right)&\leq& \II\left(|\hat{\beta}_{\omega} -\beta_{\omega} |> \frac{1}{2}\lambda^j_{\eps,\delta} \right) + \II\left(|\beta_{\omega} |>\frac{1}{2}\lambda^j_{\eps,\delta}\right).
\eeqn
Thus, the remaining two terms  $\mathbb{E}_2$ and $\mathbb{E}_3$ can be partitioned as follows
\be  \label{R_2}
\mathbb{E}_2\leq \mathbb{E}_{21} +\mathbb{E}_{22}+ \mathbb{E}_{23}, \  \mathbb{E}_3 \leq \mathbb{E}_{31}+\mathbb{E}_{32} + \mathbb{E}_{33}+ \mathbb{E}_{34},
\ee  where
\beqn  
\mathbb{E}_{21}&=&\sum\limits_{\omega \in \Omega(J,J')} \mathbb{E}\left[|\hat{\beta}_{\omega}- \beta_{\omega} |^2\II(\Omega_1)\II(\Omega_2)\II(\Theta_{\omega,\frac{1}{2}})\right],\\
\mathbb{E}_{22}&=&\sum\limits_{\omega \in \Omega(J,J')} \mathbb{E}\left[|\hat{\beta}_{\omega}- \beta_{\omega} |^2\II(\Omega_1)\II(\Omega_2)\II\left(|\beta_{\omega} |> \frac{1}{2}\lambda^j_{\eps,\delta}\right)\right],\\
\mathbb{E}_{23}&=&  \sum\limits_{\omega \in \Omega(J,J')}\EE \left[\left| \hat{ \beta}_{\omega}-\beta_{\omega}\right|^2\II(\Omega_1)\II(\Omega^c_2) \II \left(  \left|  \hat{ \beta}_{\omega} \right| >  \lambda^{j}_{\varepsilon, \delta} \right)  \right],\label{r22}\\
\mathbb{E}_{31}&=&\sum\limits_{\omega \in \Omega(J,J')} \mathbb{E}\left[| \beta_{\omega}|^2 \II(\Omega_1)\II(\Omega_2)\II(\Theta_{\omega,\frac{1}{2}}) \right],\\
\mathbb{E}_{32}&=&\sum\limits_{\omega \in \Omega(J,J')} | \beta_{\omega}|^2 \II(\Omega_1)\II(\Omega_2)\II\left(|\beta_{\omega} |<\frac{3}{2}\lambda^j_{\eps,\delta}\right),\\
\mathbb{E}_{33}&=& \sum\limits_{\omega \in \Omega(J,J')}\EE \left[ \left| { \beta}_{\omega}  \right|^2\II(\Omega^c_1)\II \left(  \left|  \hat{\beta}_{\omega} \right| <  \lambda^j_{\varepsilon, \delta} \right)\right],\label{r32}\\
\mathbb{E}_{34}&=& \sum\limits_{\omega \in \Omega(J,J')}\EE \left[ \left| { \beta}_{\omega}  \right|^2\II{(\Omega^c_2)}\II \left(  \left|  \hat{\beta}_{\omega} \right| <  \lambda^j_{\varepsilon, \delta} \right)\right].\label{r33}\eeqn
Now, choose $\gamma_1$ and $\gamma_2$ in \fr{tau} so that $\tau_1,\tau_2\geq 5$. and apply Cauchy-Schwarz inequality, with \fr{prob}, $\alpha=\frac{1}{2}$, \fr{e2} and \fr{e4}, one has
\beqn
\mathbb{E}_{21}+\mathbb{E}_{31}
&\leq&\sum\limits_{\omega \in \Omega(J,J')}\sqrt{\mathbb{E}|\hat{\beta}_{\omega}- \beta_{\omega} |^4\Pr(\Theta_{\omega,\frac{1}{2}})}+\sum\limits_{\omega \in \Omega(J,J')} | \beta_{\omega}|^2 \Pr(\Theta_{\omega,\frac{1}{2}})  \nonumber\\
&\leq&CA^22^{J(2\nu+1)} 2^{\frac{3J'}{2}}  \max\{\eps^2,\delta^2\}\max \{ \eps^{\tau_1}, \delta^{ \tau_2}\} =O\left( A^2\chi^d_{\eps,\delta, A}\right).
\eeqn
For $\mathbb{E}_{33}$, applying Cauchy-Schwarz inequality and \fr{pc1}, yields
\beqn \label{e33}
\mathbb{E}_{33}
&\leq&\sum\limits_{\omega \in \Omega(J,J')} \beta^2_{\omega}  \sqrt{\Pr(\Omega^c_1) \Pr\left(|\hat{\beta}_{\omega} |<\lambda^j_{\eps,\delta}\right) } \leq\sum\limits_{\omega \in \Omega(J,J')}\beta^2_{\omega}  \sqrt{\Pr(\Omega^c_1) } \nonumber\\
&=&O\left( A^2 \delta^2 \right)=O\left( A^2\chi^d_{\eps,\delta, A}\right).
\eeqn
Now, combining $\mathbb{E}_{23}$ and $\mathbb{E}_{34}$, applying Cauchy-Schwarz inequality and choose $\rho^2\kappa^2\geq 20$, yields
\beqn
\mathbb{E}_{23}+\mathbb{E}_{34}
&=& O \left(  \sum\limits_{\omega \in \Omega(J,J')}  \sqrt{\EE \left[\left| \hat{ \beta}_{\omega}-\beta_{\omega}\right|^4\II{(\Omega_1)}\right] \Pr \left(\Omega^c_2\right)}+\sum\limits_{\omega \in \Omega(J,J')}\left| { \beta}_{\omega}  \right|^2\Pr \left(\Omega^c_2\right) \right) \nonumber\\
&=& O\left( A^2\chi^d_{\eps,\delta, A}\right).
\eeqn
For the sum of $\mathbb{E}_{22}$ and $ \mathbb{E}_{32}$, using \fr{e2} and \fr{lambda}, yields
\be
\Delta=\mathbb{E}_{22}+\mathbb{E}_{32}=O \left( \sum\limits_{j=m_0}^{J-1} \sum\limits_{j'=m_0'}^{J'-1}  \sum\limits_{k,k'} \min \left\{   \beta^2_{j,k,j',k'} , 2^{2j\nu }\max\{ \eps^2\ln(1/\eps), \delta^2 \ln^2 \delta\}     \right\} \right).
\ee
The result of the case when $\eps^2\ln(1/\eps)\geq \delta^2 \ln^2 \delta$ has been derived by Benhaddou et al.~(2013), so we skip it. It remains to study the case when
$\max\left\{ \eps^2\ln(1/\eps), \delta^2 \ln^2 \delta \right\}=\delta^2 \ln^2 \delta$.
Then $\Delta$ can be partitioned as $\Delta \leq \Delta_1+\Delta_2+\Delta_3$, where 
\beqn
\Delta_1&=&O \left( \left\{ \sum\limits_{j=j_0+1}^{J-1} \sum\limits_{j'=m_0'}^{J'-1}+\sum\limits_{j=m_0}^{J-1} \sum\limits_{j'=j_0'+1}^{J'-1} \right\}   A^2 2^{-2js'_1-2j's'_2}  \right),\\
\Delta_2&=&O \left(  \sum\limits_{j=m_0}^{j_0} \sum\limits_{j'=m_0'}^{j'_0}  2^{j(2\nu+1)+j'}\delta^2 \ln^2\delta \  \II\left( 2^{j(2\nu+1)+j'}   \leq \chi^{d-1}_{\eps,\delta, A} \right)  \right),\\
\Delta_3&=&O \left(\sum\limits_{j=m_0}^{j_0} \sum\limits_{j'=m_0'}^{j'_0}   A^{p'} 2^{-p'js'_1-p'j's'_2} ( 2^{2j\nu }\delta^2 \ln^2 \delta  )^{1-p'/2} \II\left( 2^{j(2\nu+1)+j'}  > \chi^{d-1}_{\eps,\delta, A}    \right) \right).
\eeqn
Combining $\Delta_1$ and $\Delta_2$, and keeping in mind $j_0$ and $j'_0$ given by \fr{chijj}, one has 
\be
\Delta_1+ \Delta_2=O\left( A^2\chi^d_{\eps,\delta, A}\right).
\ee
For $\Delta_3$, we need to consider three different cases.\\
\underline{Case 1: $s_2(2\nu+1)\leq s_1$}. Then $d=\frac{2s_2}{2s_2+1}$, and 
\be
\Delta_3=O\left( A^2\left( \chi_{\eps, \delta, A} \right)^{\frac{2s_2}{2s_2+1}} \sum\limits_{j=m_0}^{j_0}2^{ -jp' \left[s_1- (2\nu+1)s_2 \right]} \right)=O \left(A^2\chi^d_{\eps,\delta, A}   [\ln (1/\delta)]^{\II( s_1=s_2(2\nu+1))}       \right).
\ee
\underline{Case 2: $  (2\nu+1)\left(\frac{1}{p}-\frac{1}{2}\right)< s_1 < s_2(2\nu+1) $}. Then $d=\frac{2s_1}{2s_1+2\nu+1}$, and 
\be
\Delta_3=O\left( A^2 \left(\chi_{\eps, \delta, A}\right)^{\frac{2s_1}{2s_1+2\nu+1}} \sum\limits_{j'=m_0'}^{j'_0} 2^{-j'p'\left(s_2-\frac{s_1}{2\nu+1}\right)} \right)=O \left( A^2 \chi^d_{\eps, \delta, A} \right).
\ee
\underline{Case 3: $s_1 \leq (2\nu+1)\left(\frac{1}{p}-\frac{1}{2}\right)$}. Then $d=\frac{2s'_1}{2s'_1+2\nu}$, and 
\be \label{del33}
\Delta_3=\left(  A^2 \left( \chi_{\eps,\delta, A} \right)^{1-\frac{p'}{2}}\sum\limits_{j=m_0}^{j_0}  2^{-j\left(ps'_1-2\nu\left(1-\frac{p}{2} \right)   \right)}\right)=O \left( A^2 \chi^d_{\eps, \delta, A}    [\ln (1/\delta)]^{\II\left(s_1 =(2\nu+1)\left(\frac{1}{p}-\frac{1}{2}\right) \right)}    \right).
\ee
Combining \fr{e1234}-\fr{del33} completes the proof. $\Box$
\bibliographystyle{model1-num-names}

\begin{thebibliography}{}
\bibitem{abr}Abramovich, F., Silverman, B.W., Wavelet decomposition approaches to statistical inverse problems. Vol. {\it Biometrika}, {\bf85}(1),  115--129, 1998.

 \bibitem{ben 1} Benhaddou, R., Minimax lower bounds for the simultaneous wavelet deconvolution with fractional Gaussian noise and unknown kernels. {\it Statistics and Probability Letters}, Vol. {\bf 140} (C), 91-95, 2018b.

\bibitem{ben2}
Benhaddou, R., Blind deconvolution model in periodic setting with fractional Gaussian noise. \textit{Communication in Statistics-Theory and Methods}, Vol. \textbf{0}, 1-12, 2018a.

\bibitem{ben3}
Benhaddou, R.,
On minimax convergence rates under $L^p$-risk for the anisotropic  functional deconvolution model.
 {\it Statistics and Probability Letters,} Vol. {\bf 130},  120-125, 2017.
 
\bibitem{ben4}
Benhaddou, R., Pensky, M. and Picard, D., Anisotropic de-noising in functional deconvolution model with dimension-free convergence rates. \textit{Electronic Journal of Statistics}, Vol. \textbf{7}, 1686-1715, 2013.

\bibitem{ben5}
Benhaddou, R., Pensky, M. and Rajapakshage, R., Anisotropic functional Laplace deconvolution. \textit{Journal of Statistical Planning and Inference}, Vol. \textbf{199}, 271-285, 2019.
\bibitem{}
Bunea, F., Tsybakov, A. and Wegkamp, M. H., Aggregation for Gaussian regression. \textit{Ann. Statist.}, Vol. \textbf{35}, 1674-1697, 2007.

\bibitem{}
Delattre, S., Hoffmann, M., Picard, D., Vareschi, T., Blockwise SVD with error in the operator and application to blind deconvolution. \textit{Electronic Journal of Statistics}, Vol. \textbf{6}, 2274-2308, 2012.

\bibitem{}
Donoho, D. L., Nonlinear solution of linear inverse problems by wavelet-vaguelette decomposition. \textit{Appl. Computat. Harmon. Anal.}, Vol. \textbf{2}, 101-126, 1995. 

\bibitem{}
Donoho, D. L. and Johnstone, I. M., Ideal Spatial Adaptation by Wavelet Shrinkage. \textit{Biometrika} \textbf{81} (3), 425-455, 1994.



\bibitem  {donrai}
Donoho, D.L., Raimondo, M., Translation invariant
deconvolution in a periodic setting. {\it International Journal of
Wavelets, Multiresolution and Information Processing}, Vol. {\bf 14},
415--432, 2004.


\bibitem{dozzi}
Dozzi, M., {\it Stochastic Processes with a Multidimensional Parameter,} Longman
Scientific \& Technical, New York, 1989.

\bibitem{}
Hoffmann, M., Reiss, M., Nonlinear estimation for linear inverse problems with error in the operator.  \textit{Ann. Statist.},Vol. \textbf{36}, 310-336, 2008. 

\bibitem{}
Johnston, I.M., Kerkyacharian, G., Picard, D. and Raimondo, M., Wavelet Deconvolution in a Periodic Setting. \textit{Journal of the Royal Statistical Society. Series B (Statistical Methodology)}, Vol. \textbf{66} (3), 547-573, 2004.
\bibitem{}
\O ksendal, B., \textit{Stochastic Differential Equations: An Introduction with Applications}, 6th edition. Springer, 2003.


\bibitem{}
Pensky, M., Sapatinas, T., Functional deconvolution in a periodic setting: uniform case. \textit{Ann. Statist.}, Vol. \textbf{37}, 73-104, 2009.

\bibitem{}
Pensky, M., Vidakovic, B., Adaptive wavelet estimator for nonparametric densify deconvolution. \textit{Ann. Statist.}, Vol. \textbf{37}, 2033-2053, 1999.

\bibitem{}
Robinson, E. A., \textit{Seismic Inversion and Deconvolution: Part B: Dual-Sensor Technology}, Elsevier, Oxford, 1999.

\bibitem{}
Tsybakov, A.B., \textit{Introduction to Nonparametric Esitmation}. Springer, New York, 2009.

\bibitem{}
Vareschi, T., Noisy Laplace deconvolution with error in the operator. \textit{Journal of Statistical Planning and Inference}, Vol. \textbf{157-158}, 16-35, 2015.

\bibitem {walter}
Walter, G., Shen, X., Deconvolution using Meyer wavelets.
{\it Journal of Integral Equations and Applications}, Vol. {\bf 11},
515--534, 1999.

\end{thebibliography}

\end{document}